\documentclass[12pt,reqno]{amsart}
\usepackage{amsfonts,amsmath,amssymb}

\allowdisplaybreaks
\advance\textheight by \topskip

\usepackage{mathrsfs}

\usepackage[latin1]{inputenc}
\usepackage{graphicx}

\usepackage{tikz}
\newtheorem{theorem}{Theorem}
\usepackage{euscript}

\def\[{[\! [}
\def\]{]\! ]}

\title[Unary binary trees, Hex-trees, and marked ordered trees]{Weighted unary-binary trees, Hex-trees, marked ordered trees, and related structures}

\author[H.~Prodinger]{Helmut Prodinger}

\address{Helmut Prodinger,
Mathematics Department, Stellenbosch University,
7602 Stellenbosch, South Africa.}
\email{hproding@sun.ac.za}

\thanks{Results from https://arxiv.org/abs/2105.03350 are now included in the present version (August 16, 2021).}

\date{\today}

\begin{document}

\begin{abstract}
	Hex-trees are identified as a particular instance of weighted unary-binary trees. The Horton-Strahler 
	numbers of these objects are revisited, and, thanks to a substitution that is not immediately intuitive, explicit
	results are possible. They are augmented by asymptotic evaluations as well. Furthermore, marked ordered trees
	(in bijection to skew Dyck paths) are investigated, followed 3-Motzkin paths and multi-edge trees. The underlying
	theme is sequence A002212 in the Encyclopedia of integer sequences
\end{abstract}

\maketitle

\section{Introduction}

This paper gives some (mostly new) results about the sequence
\begin{equation*}
1, 1, 3, 10, 36, 137, 543, 2219, 9285, 39587, 171369, 751236, 3328218, 14878455,\dots,
\end{equation*}
which is A002212 in \cite{OEIS}.

Here is the plan about the structures enumerated by this sequence:

Hex-trees \cite{KimStanley}; they are identified as weighted unary-binary trees, with weight one. Apart from left and right branches, as
in binary trees, there are also unary branches, and they can come in different colours, here in just one colour.
Unary-binary trees played a role in the present authors scientific development, as documented in \cite{FlPr86}, a paper written with
the late and great Philippe Flajolet, about the register function (Horton-Strahler numbers) of unary-binary trees. Here, we can offer an improvement, using
a ``better'' substitution than in \cite{FlPr86}. The results can now be made fully explicit. As a by-product, this provides a definition and analysis
of the Horton-Strahler numbers of Hex-trees. An introductory section (about binary trees) provides all the basics.

Then we move to skew Dyck paths \cite{Deutsch-italy}. They are like Dyck paths, but allow for an extra step $(-1,-1)$, provided that the path does not intersect itself.
An equivalent model, defined and described using a bijection, is from \cite{Deutsch-italy}: Marked ordered trees. They are like ordered trees, with an additional feature,
namely each rightmost edge (except for one that leads to a leaf) can be coloured with two colours. Since we find this class of trees to be interesting, we
analyze two parameters of them: number of leaves and height. While the number of leaves for ordered trees is about $n/2$, it is only $n/10$ in the new model.
For the height, the leading term $\sqrt{\pi n}$ drops to $\frac{2}{\sqrt 5}\sqrt{\pi n}$. Of course, many more parameters of this new class of trees could be investigated,
which we encourage to do.

The last two classes of structures are multi-edge trees. Our interest in them was already triggered in an earlier publication \cite{HPW}. They may be seen as ordered trees, but
with weighted edges. The weights are integers $\ge1$, and a weight $a$ may be interpreted as $a$ parallel edges. The other class are 3-Motzkin paths. They are like
Motzkin paths (Dyck paths plus horizontal steps); however, the horizontal steps come in three different colours. A bijection is described. Since 3-Motzkin paths and multi-edge trees
 are very much alike (using a variation of the classical rotation correspondence), all the structures that are discussed in this paper can be linked via bijections.

\section{Binary trees and Horton-Strahler numbers}

Binary trees may be expressed by the following symbolic equation, which says that they include the empty tree 
and trees recursively built from a root followed by two subtrees, which are binary trees:
\begin{center}\small
	\begin{tikzpicture}
	[inner sep=1.3mm,
	s1/.style={circle=10pt,draw=black!90,thick},
	s2/.style={rectangle,draw=black!50,thick},scale=0.5]
	
	\node at ( -4.8,0) { $\mathscr{B}$};
	
	\node at (-3,0) { $=$};
	\node(c) at (-1.5,0){ $\qed$};
	\node at (0.7,0) {$+$};
	\node(d) at (3,1)[s1]{};
	\node(e) at (2,-1){ $\mathscr{B}$};
	\node(f) at (4,-1){ $\mathscr{B}$};
	\path [draw,-,black!90] (d) -- (e) node{};%
	\path [draw,-,black!90] (d) -- (f) node{};%
	
	\end{tikzpicture}
\end{center}

Binary trees are counted by Catalan numbers and there is an important  parameter \textsf{reg}, which 
in Computer Science is called the register function. It associates to each binary tree (which is used to code an arithmetic expression,
with data in the leaves and operators in the internal nodes) the minimal number of extra registers that is needed to evaluate the tree. The optimal strategy is to evaluate difficult subtrees first, and use one register to keep its value, which does not hurt, if the other subtree requires less registers. If both subtrees are equally difficult, then one more register is used, compared to the requirements of the subtrees. This natural parameter is among Combinatorialists known as the Horton-Strahler numbers, and we will adopt this
name throughout this paper.

There is a recursive description of this function: $\textsf{reg}(\square)=0$, and if tree $t$ has subtrees $t_1$ and $t_2$, then
\begin{equation*}
\textsf{reg}(t)=
\begin{cases}
\max\{\textsf{reg}(t_1),\textsf{reg}(t_2)\}&\text{ if } \textsf{reg}(t_1)\ne\textsf{reg}(t_2),\\
1+\textsf{reg}(t_1)&\text{ otherwise}.
\end{cases}
\end{equation*}

The recursive description attaches numbers to the nodes, starting with 0's at the leaves and
then going up; the number appearing at the root is the Horton-Strahler number of the tree.
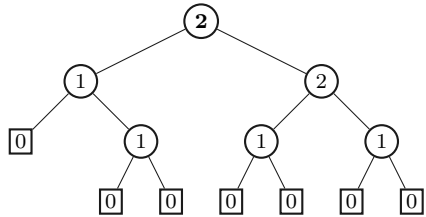
\begin{center}\tiny
\begin{tikzpicture}
[scale=0.4,inner sep=0.7mm,
s1/.style={circle,draw=black!90,thick},
s2/.style={rectangle,draw=black!90,thick}]
\node(a) at ( 0,8) [s1] [text=black]{$\boldsymbol{2}$};
\node(b) at ( -4,6) [s1] [text=black]{$1$};
\node(c) at ( 4,6) [s1] [text=black]{$2$};
\node(d) at ( -6,4) [s2] [text=black]{$0$};
\node(e) at ( -2,4) [s1] [text=black]{$1$};
\node(f) at ( 2,4) [s1] [text=black]{$1$};
\node(g) at ( 6,4) [s1] [text=black]{$1$};
\node(h) at ( -3,2) [s2] [text=black]{$0$};
\node(i) at ( -1,2) [s2] [text=black]{$0$};
\node(j) at ( 1,2) [s2] [text=black]{$0$};
\node(k) at ( 3,2) [s2] [text=black]{$0$};
\node(l) at ( 5,2) [s2] [text=black]{$0$};
\node(m) at ( 7,2) [s2] [text=black]{$0$};
\path [draw,-,black!90] (a) -- (b) node{};%
\path [draw,-,black!90] (a) -- (c) node{};%
\path [draw,-,black!90] (b) -- (d) node{};%
\path [draw,-,black!90] (b) -- (e) node{};%
\path [draw,-,black!90] (c) -- (f) node{};%
\path [draw,-,black!90] (c) -- (g) node{};%
\path [draw,-,black!90] (e) -- (h) node{};%
\path [draw,-,black!90] (e) -- (i) node{};%
\path [draw,-,black!90] (f) -- (j) node{};%
\path [draw,-,black!90] (f) -- (k) node{};%
\path [draw,-,black!90] (g) -- (l) node{};%
\path [draw,-,black!90] (g) -- (m) node{};%
\end{tikzpicture}
\end{center}

Let $\mathscr{R}_{p}$ denote the family of
trees with Horton-Strahler number $=p$, then one gets immediately from the recursive 
definition:
\begin{center}\small
\begin{tikzpicture}
[inner sep=1.3mm,
s1/.style={circle=10pt,draw=black!90,thick},
s2/.style={rectangle,draw=black!50,thick},scale=0.5]

\node at ( -5,0) { $\mathscr{R}_p$};

\node at (-4,0) { $=$};
\node(a) at (-2,1)[s1]{};
\node(b) at (-3,-1){ $\mathscr{R}_{p-1}$};
\node(c) at (-1,-1){ $\mathscr{R}_{p-1}$};
\path [draw,-,black!90] (a) -- (b) node{};%
\path [draw,-,black!90] (a) -- (c) node{};%
\node at (0.7,0) {$+$};
\node(d) at (3,1)[s1]{};
\node(e) at (2,-1){ $\mathscr{R}_{p}$};
\node(f) at (4,-1.2){ $\sum\limits_{j<p}\mathscr{R}_{j} $};
\path [draw,-,black!90] (d) -- (e) node{};%
\path [draw,-,black!90] (d) -- (f) node{};%
\node at (5+0.7,0) {$+$};
\node(dd) at (5.5+3,1)[s1]{};
\node(ee) at (5.5+2,-1.2){ $\sum\limits_{j<p}\mathscr{R}_{j}$};
\node(ff) at (5.5+4,-1){ $\mathscr{R}_{p}$};
\path [draw,-,black!90] (dd) -- (ee) node{};%
\path [draw,-,black!90] (dd) -- (ff) node{};%
\end{tikzpicture}
\end{center}
 In terms of generating functions, these equations read as
\begin{equation*}
R_p(z)=zR_{p-1}^2(z)+2zR_p(z)\sum_{j<p}R_j(z);
\end{equation*}
the variable $z$ is used to mark the size (i.~e., the number of internal nodes) of the binary tree.

A historic account of these concepts, from the angle of Philippe Flajolet, who was one of the pioneers
is \cite{register-introduction}, compare also \cite{ECA-historic}.

Amazingly, the recursion for the generating functions $R_p(z)$  can be solved explicitly! The substitution
\begin{equation*}
z=\frac{u}{(1+u)^2}
\end{equation*}
that de Bruijn, Knuth, and Rice~\cite{BrKnRi72} also used, produces the nice expression
\begin{equation*}
R_p(z)=\frac{1-u^2}{u}\frac{u^{2^p}}{1-u^{2^{p+1}}}.
\end{equation*}
Of course, once this is \emph{known}, it can be proved by induction, using the recursive formula. For the readers benefit, this will be sketched now. 

We start with the auxiliary formula
\begin{equation*}
\sum_{0\le j<p}\frac{u^{2^j}}{1-u^{2^{j+1}}}=\frac{u}{1-u}-\frac{u^{2^p}}{1-u^{2^{p}}},
\end{equation*}
which is easy to prove by induction: For $p=0$, the formula $0=\frac{u}{1-u}-\frac{u}{1-u}$ is correct, and then
\begin{align*}
\sum_{0\le j<p+1}
\frac{u^{2^j}}{1-u^{2^{j+1}}}
	&=\frac{u}{1-u}-\frac{u^{2^p}}{1-u^{2^{p}}}+\frac{u^{2^p}}{1-u^{2^{p+1}}}\\
	&=\frac{u}{1-u}-\frac{u^{2^p}(1+u^{2^p})}{1-u^{2^{p+1}}}+\frac{u^{2^p}}{1-u^{2^{p+1}}}
	=\frac{u}{1-u}-\frac{u^{2^{p+1}}}{1-u^{2^{p+1}}}.
		\end{align*}
Now the formula for $R_p(z)$ can also be proved by induction. First, $R_0(z)=\frac{1-u^2}{u}\frac{u}{1-u^{2}}=1$, as it should, and
\begin{align*}
zR_{p-1}^2(z)&+2zR_p(z)\sum_{j<p}R_j(z)\\
&=\frac{u}{(1+u)^2}\frac{(1-u^2)^2}{u^2}\frac{u^{2^{p}}}{(1-u^{2^{p}})^2}
+\frac{2u}{(1+u)^2}R_p(z)\sum_{j<p}\frac{1-u^2}{u}\frac{u^{2^j}}{1-u^{2^{j+1}}}\\
&= \frac{u^{2^{p}-1}(1-u)^2}{(1-u^{2^{p}})^2}
+\frac{2(1-u)}{(1+u)}R_p(z)\sum_{j<p}\frac{u^{2^j}}{1-u^{2^{j+1}}}.
\end{align*}
Solving
\begin{align*}
R_p(z)= \frac{u^{2^{p}-1}(1-u)^2}{(1-u^{2^{p}})^2}
+\frac{2(1-u)}{(1+u)}R_p(z)\bigg[\frac{u}{1-u}-\frac{u^{2^p}}{1-u^{2^{p}}}\bigg]
\end{align*}
leads to
\begin{align*}
R_p(z)\frac{1-u}{1+u}\bigg[1+2\frac{u^{2^p}}{1-u^{2^{p}}}\bigg]= \frac{u^{2^{p}-1}(1-u)^2}{(1-u^{2^{p}})^2},
\end{align*}
or, simplified
\begin{align*}
R_p(z)= \frac{u^{2^{p}-1}(1-u^2)}{(1-u^{2^{p}})(1+u^{2^{p}})}
=\frac{1-u^2}{u}\frac{u^{2^{p}}}{(1-u^{2^{p+1}})},
\end{align*}
which is the formula that we needed to prove. \qed

 \section{Unary-binary trees and Hex-trees}
The family of unary-binary trees ${\mathscr{M}}$ might be defined by the symbolic equation
\begin{center}
	\begin{tikzpicture}
	[inner sep=1.3mm,
	s1/.style={circle=10pt,draw=black!90,thick},
	s2/.style={rectangle,draw=black!50,thick},scale=0.5]
	
	\node at ( -12.8,0.1) { ${\mathscr{M}}$};
	
	\node at (-11.2,0) { $=$};
	\node at (-4.5,0) { $+$};
	\node(a) at (-2,1)[s1]{};
	\node(b) at (-3,-1){ ${\mathscr{M}}$};
	\node(c) at (-1,-1){ ${\mathscr{M}}$};
	\path [draw,-,black!90] (a) -- (b) node{};%
	\path [draw,-,black!90] (a) -- (c) node{};%
	\node(a1) at (-6.5,1)[s1]{};
	\node(b1) at (-6.5,-1){ $ {\mathscr{M}}\setminus\{\square\}$};
	\path [draw,-,black!90] (a1) -- (b1) node{};%
	\node at (-9.0,0) { $\square\ \ +$};
	\end{tikzpicture}
\end{center}
The equation for the generating function is
\begin{equation*}
M=1+z(M-1)+zM^2
\end{equation*}
 with the solution
 \begin{equation*}
M=M(z)=\frac{1-z-\sqrt{1-6z+5z^2}}{2z}=1+z+3{z}^{2}+10{z}^{3}+36{z}^{4}+\cdots;
 \end{equation*}
the coefficients form again sequence A002212 in \cite{OEIS} and enumerate Schr\"oder paths, among many other things.
We will come to equivalent structures a bit later. 

In the instance of unary-binary trees, we can also work with a substitution: Set $z=\frac{u}{1+3u+u^2}$, then
$M(z)=1+u$. Unary-binary trees and the register function were investigated in \cite{FlPr86}, but the present favourable substitution was not used.
Therefore, in this previous paper, asymptotic results were available but no explicit formulae.

This works also with a weighted version, where we allow unary edges with $a$ different colours. Then
\begin{center}
	\begin{tikzpicture}
	[inner sep=1.3mm,
	s1/.style={circle=10pt,draw=black!90,thick},
	s2/.style={rectangle,draw=black!50,thick},scale=0.5]
	
	\node at ( -12.8,0.1) { ${\mathscr{N}}$};
	
	\node at (-11.5,0) { $=$};
	\node at (-4.5,0) { $+$};
	\node(a) at (-2,1)[s1]{};
	\node(b) at (-3,-1){ ${\mathscr{N}}$};
	\node(c) at (-1,-1){ ${\mathscr{N}}$};
	\path [draw,-,black!90] (a) -- (b) node{};%
	\path [draw,-,black!90] (a) -- (c) node{};%
\node at (-7.5,0){$a\ \cdot$};
	\node(a1) at (-6.5,1)[s1]{};
	\node(b1) at (-6.5,-1){ $ {\mathscr{N}}\setminus\{\square\}$};
	\path [draw,-,black!90] (a1) -- (b1) node{};%
	\node at (-9.5,0) { $\square\ \ +$};
	\end{tikzpicture}
\end{center}
and with the substitution $z=\frac{u}{1+(a+2)u+u^2}$, the generating function is beautifully expressed as $N(z)=1+u$.
For $a=0$, this covers also binary trees.

We will consider  the Horton-Strahler numbers of unary-binary trees in the sequel. The definition is naturally extended by

\begin{center}
\vspace{0pt}\begin{tikzpicture}
	[inner sep=0.6mm,
s1/.style={circle=1pt,draw=black!90,thick}]
\node[] at ( -0.300,-0.10) {\textsf{reg}\bigg(};
\node[] at ( 0.6500,-0.10) {\bigg)};
\node[] at ( 1.5500,-0.10) {=\ \textsf{reg}(t).};
\path [draw,-,black!90 ] (0.3,0.34) -- (0.3,-0.350) ;
\node [s1]at ( 0.300,0.4) { };
\node[] at ( 0.300,-0.60) {$t$ };
 \end{tikzpicture} 
\end{center}

Now we can move again to $R_p(z)$, the generating funciton of (generalized) unary-binary trees with Horton-Strahler number $=p$.
The recursion (for $p\ge1$) is
\begin{center}\small
	\begin{tikzpicture}
	[inner sep=1.3mm,
	s1/.style={circle=10pt,draw=black!90,thick},
	s2/.style={rectangle,draw=black!50,thick},scale=0.5]
	
	\node at ( -5,0) { $\mathscr{R}_p$};
	
	\node at (-4,0) { $=$};
	\node(a) at (-2,1)[s1]{};
	\node(b) at (-3,-1){ $\mathscr{R}_{p-1}$};
	\node(c) at (-1,-1){ $\mathscr{R}_{p-1}$};
	\path [draw,-,black!90] (a) -- (b) node{};%
	\path [draw,-,black!90] (a) -- (c) node{};%
	\node at (0.7,0) {$+$};
	\node(d) at (3,1)[s1]{};
	\node(e) at (2,-1){ $\mathscr{R}_{p}$};
	\node(f) at (4,-1.2){ $\sum\limits_{j<p}\mathscr{R}_{j} $};
	\path [draw,-,black!90] (d) -- (e) node{};%
	\path [draw,-,black!90] (d) -- (f) node{};%
	\node at (5+0.7,0) {$+$};
	\node(dd) at (5.5+3,1)[s1]{};
	\node(ee) at (5.5+2,-1.2){ $\sum\limits_{j<p}\mathscr{R}_{j}$};
	\node(ff) at (5.5+4,-1){ $\mathscr{R}_{p}$};
	\path [draw,-,black!90] (dd) -- (ee) node{};%
	\path [draw,-,black!90] (dd) -- (ff) node{};%
	\node(dd) at (13,1)[s1]{};
	\node(ee) at (13,-1){ $\mathscr{R}_{p}$};
	\path [draw,-,black!90] (dd) -- (ee) node{};%
	\node at (11.5,0) {$+\ \ a\cdot$};
	\end{tikzpicture}
\end{center}
In terms of generating functions, these equations read as
\begin{equation*}
R_p(z)=zR_{p-1}^2(z)+2zR_p(z)\sum_{j<p}R_j(z)+azR_p(z), \quad p\ge1;\quad R_0(z)=1.
\end{equation*}
Amazingly, with the substitution $z=\frac{u}{1+(a+2)u+u^2}$, formally we get the \emph{same} solution as in the binary case:
\begin{equation*}
R_p(z)=\frac{1-u^2}{u}\frac{u^{2^p}}{1-u^{2^{p+1}}}.
\end{equation*}
The proof by induction is as before. One sees another advantage of the substitution: On a formal level, many manipulations do not need to be repeated. Only when one switches back to the $z$-world, things become different.

Now we move to Hex-trees.
\begin{center}
	\begin{tikzpicture}
	[inner sep=1.3mm,
	s1/.style={circle=10pt,draw=black!90,thick},
	s2/.style={rectangle,draw=black!50,thick},scale=0.5]
	
	\node at ( -12.4,0.1) { ${\mathscr{H}}$};
	
	\node at (-11.0,0) { $=$};
	
	\node(a) at (-1,1)[s1]{};
	\node(b) at (-3,-1){ ${\mathscr{H}\setminus\{\square\}}$};
	\node(c) at (1,-1){ ${\mathscr{H}\setminus\{\square\}}$};
	\path [draw,-,black!90] (a) -- (b) node{};%
	\path [draw,-,black!90] (a) -- (c) node{};%

\begin{scope}[xshift=13cm]
\node at (-10.0,0) { $+$};
\node at (-4.5,0) { $+$};
\node(a1) at (-6.5,1)[s1]{};
\node(b1) at (-7.5,-1){ $ {\mathscr{H}}\setminus\{\square\}$};
\path [draw,-,black!90] (a1) -- (b1) node{};%
\end{scope}

\begin{scope}[xshift=21cm]
\node(a1) at (-6.5,1)[s1]{};
\node(b1) at (-5.5,-1){ $ {\mathscr{H}}\setminus\{\square\}$};
\path [draw,-,black!90] (a1) -- (b1) node{};%
\end{scope}

	\begin{scope}[xshift=17cm]
	\node at (-4.5,0) { $+$};
	\node(a1) at (-6.5,1)[s1]{};
	\node(b1) at (-6.5,-1){ $ {\mathscr{H}}\setminus\{\square\}$};
	\path [draw,-,black!90] (a1) -- (b1) node{};%
	\end{scope}

	\node at (-9.0,0) { $\square\ \ +$};
	\node[s1] at (-7.0,0) { };
	\node at (-5.0,0) {$+$ };
	\end{tikzpicture}
\end{center}

Hex trees either have two non-empty successors, or one of 3 types of unary successors (called left, middle, right).
The author has seen this family first in \cite{KimStanley}, but one can find older literature following the references and the usual
search engines.

The generating function satisfies
\begin{align*}
H&(z)=1+z(H(z)-1)^2+z+3z(H(z)-1)=\frac{1-z-\sqrt{(1-z)(1-5z)}}{2z}\\
&=1+z+3{z}^{2}+10{z}^{3}+36{z}^{4}+137{z}^{5}+543{z}^{6}+2219
{z}^{7}+9285{z}^{8}+39587{z}^{9}+\cdots.
\end{align*}
The same generating function also appears in \cite{HPW}, and it is again sequence A002212 in \cite{OEIS}.
One can rewrite the symbolic equation as
\begin{center}
	\begin{tikzpicture}
	[inner sep=1.3mm,
	s1/.style={circle=10pt,draw=black!90,thick},
	s2/.style={rectangle,draw=black!50,thick},scale=0.5]
	
	\node at ( -12.4,0.1) { ${\mathscr{H}}$};
	
	\node at (-11.0,0) { $=$};
	
		\begin{scope}[xshift=-4cm]
	\node(a) at (-1,1)[s1]{};
	\node(b) at (-3,-1){ $\mathscr{H}$};
	\node(c) at (1,-1){ $\mathscr{H}$};
	\path [draw,-,black!90] (a) -- (b) node{};%
	\path [draw,-,black!90] (a) -- (c) node{};%
		\end{scope}

	\begin{scope}[xshift=7.5cm]
	\node at (-8.5,0) { $+$};
	\node(a1) at (-6.5,1)[s1]{};
	\node(b1) at (-6.5,-1){ $ {\mathscr{H}}\setminus\{\square\}$};
	\path [draw,-,black!90] (a1) -- (b1) node{};%
	\end{scope}

	\node at (-9.0,0) { $\square\ \ \ +$};
	
	\end{tikzpicture}
\end{center}
and sees in this way that the Hex-trees are unary-binary trees (with parameter $a=1$).

\subsection*{Continuing with enumerations}

First, we will enumerate the number of (generalized) unary-binary trees with $n$ (internal) nodes. For that we need the notion 
of generalized trinomial coefficients, viz.
\begin{equation*}
\binom{n;1,a,1}{k}:=[z^k](1+az+z^2)^n.
\end{equation*}
Of course, for $a=2$, this simplifies to a binomial coefficient $\binom{2n}{k}$. We will use contour integration to pull out coefficients, and the contour of integer, in whatever variable, is a small circle (or equivalent) around the origin.
The desired number is
\begin{align*}
[z^n](1+u)&=\frac1{2\pi i}\oint \frac{dz}{z^{n+1}}(1+u)\\
&=\frac1{2\pi i}\oint \frac{du(1-u^2)(1+(a+2)u+u^2)^{n+1}}{(1+(a+2)u+u^2)^2u^{n+1}}(1+u)\\
&=[u^{n+1}](1-u)(1+u)^2(1+(a+2)u+u^2)^{n-1}\\
&=\binom{n-1;1,a+2,1}{n+1}+\binom{n-1;1,a+2,1}{n}\\*
&\hspace*{4cm}-\binom{n-1;1,a+2,1}{n-1}-\binom{n-1;1,a+2,1}{n-2}.
\end{align*}
Then we introducte $S_p(z)=R_{p}(z)+R_{p+1}(z)+R_{p+2}(z)+\cdots$, the generating function of trees with Horton-Strahler number
$\ge p$. Using the summation formula proved earlier, we get
\begin{equation*}
S_p(z)=\frac{1-u^2}{u}\frac{u^{2^p}}{1-u^{2^{p}}}=
\frac{1-u^2}{u}\sum_{k\ge1}u^{k2^p}.
\end{equation*}
Further,
\begin{align*}
[z^n]S_p(z)&=\sum_{k\ge1}\frac1{2\pi i}\oint \frac{dz}{z^{n+1}}\frac{1-u^2}{u}u^{k2^p}.
\end{align*}

\subsection*{Asymptotics}

We start by deriving asymptotics for the number of (generalized) unary-binary trees with $n$ (internal) nodes. This is
a standard application of singularity analysis of generating functions, as described in \cite{FlOd90} and \cite{FS}.

We start from the generating function
\begin{equation*}
	N(z)=\frac{1-az-\sqrt{1-2(a+2)z+a(a+4)z^2}}{2z}
\end{equation*}
and determine the singularity closest to the origin, which is the value making the square root disappear:
 $z=\frac1{a+4}$.
After that, the local expansion of $N(z)$ around this singularity is determined:
\begin{equation*}
N(z) \sim	2-\sqrt{a+4}\sqrt{1-(a+4)z}.
\end{equation*}
The translation lemmas given in \cite{FlOd90} and \cite{FS} provide the asymptotics:
\begin{align*}
	[z^n]N(z)&\sim [z^n]\Big(2-\sqrt{a+4}\sqrt{1-(a+4)z}\Big)\\&
	=-\sqrt{a+4}(a+4)^n\frac{n^{-3/2}}{\Gamma(-\frac12)}=(a+4)^{n+1/2}\frac{1}{2\sqrt\pi n^{3/2}}.
\end{align*}
Just note that $a=0$ is the well-known formula for binary trees with $n$ nodes.

Now we move to the generating function for the average number of registers. Apart from normalization it is
\begin{align*}
\sum_{p\ge1}pR_p(z)&=\sum_{p\ge1}S_p(z)=\frac{1-u^2}{u}\sum_{p\ge1}\sum_{k\ge1}u^{k2^p}\\
&=\frac{1-u^2}{u}\sum_{n\ge1}v_2(n)u^n,
\end{align*}
where $v_2(n)$ is the highest exponent $k$ such $2^k$ divides $n$.

This has to be studied around $u=1$, which, upon setting $u=e^{-t}$, means around $t=0$.
Eventually, and that is the only thing that is different here, this is to be retranslated  into a singular expansion of $z$ around
its singularity, which depends on the parameter $a$. 

For the reader's convenience, we also repeat the steps that were known before.
The first factor is elementary:
\begin{equation*}
\frac{1-u^2}{u}\sim2t+{\frac {1}{3}}{t}^{3}+\cdots
\end{equation*}
For 
\begin{equation*}
\sum_{p\ge1}\sum_{k\ge1}e^{-k2^pt},
\end{equation*}
one applies the Mellin transform, with the result
\begin{equation*}
\frac{\Gamma(s)\zeta(s)}{2^s-1}.
\end{equation*}
Applying the inversion formula, one finds
\begin{equation*}
\sum_{p\ge1}\sum_{k\ge1}e^{-k2^pt}=\frac1{2\pi i}\int_{2-i\infty}^{2+i\infty}t^{-s}\frac{\Gamma(s)\zeta(s)}{2^s-1}ds.
\end{equation*}
Shifting the line of integration to the left, the residues at the poles $s=1$, $s=0$, $s=\chi_k=\frac{2k\pi i}{\log2}$, $k\neq0$ provide enough terms for our asymptotic expansion.
\begin{equation*}
\frac1{t}+{\frac {\gamma}{2\log2 }}-\frac14-	\frac {\log   \pi   }{2\log2 }+\frac {\log t }{2\log2}
+\frac1{\log2}\sum_{k\neq0}\Gamma(\chi_k)\zeta(\chi_k)t^{-\chi_k}.
\end{equation*}
Combined with the elementary factor, this leads to
\begin{equation*}
2+\Big(\frac {\gamma}{\log2 }-\frac12-\frac {\log   \pi   }{\log2 }+\frac {\log t }{\log2}\Big)t+\frac{2t}{\log2}\sum_{k\neq0}\Gamma(\chi_k)\zeta(\chi_k)t^{-\chi_k}+O(t^2\log t).
\end{equation*}
Now we want to translate into the original $z$-world. Since $z=\frac{u}{1+(a+2)u+u^2}$, $u=1$ translates into the singularity $z=\frac{1}{4+a}$. Further,
\begin{equation*}
t\sim \sqrt{4+a}\cdot \sqrt{1-z(4+a)},
\end{equation*}
let us abbreviate $A=4+a$, then for singularity analysis we must consider
\begin{align*}
&\frac {\sqrt{A}\cdot \sqrt{1-zA}\log (1-zA) }{2\log2}\\
&+	\Big(\frac {\gamma}{\log2 }-\frac12-\frac {\log   \pi   }{\log2 }+\frac{\log A}{2\log 2}\Big)\sqrt{A}\cdot \sqrt{1-zA}\\
&+\frac{2  }{\log2}\sum_{k\neq0}\Gamma(\chi_k)\zeta(\chi_k) A^{\frac{1-\chi_k}2}(1-zA)^{\frac{1-\chi_k}2}.
\end{align*}
The formula that is perhaps less known and needed here is \cite{FlOd90}
\begin{align*}
[z^n]\log(1-z)\sqrt{1-z}\sim \frac{n^{-3/2}\log n}{2\sqrt \pi}+\frac{n^{-3/2}}{2\sqrt \pi}(-2+\gamma +2\log2);
\end{align*}
furthermore we need
\begin{equation*}
[z^n](1-z)^\alpha \sim \frac{n^{-\alpha-1}}{\Gamma(-\alpha)}.
\end{equation*}
We start with the most complicated term:
\begin{align*}
\frac{[z^n]\frac {\sqrt{A}\cdot \sqrt{1-zA}\log (1-zA) }{2\log2}}{[z^n]N(z)}
&\sim \frac {\sqrt{A}}{2\log2}\frac{A^n\Big(\frac{n^{-3/2}\log n}{2\sqrt \pi}+\frac{n^{-3/2}}{2\sqrt \pi}(-2+\gamma +2\log2)\Big)}
{A^{n+1/2}\frac{1}{2\sqrt\pi n^{3/2}}}\\
&= \log_4 n+1+ \frac{\gamma }{2\log2}- \frac{1}{\log2}.
\end{align*}
The next term we consider is
\begin{align*}
\Big(\frac {\gamma}{\log2 }-\frac12-\frac {\log   \pi   }{\log2 }+\frac{\log A}{2\log 2}\Big)&\sqrt{A}\frac{[z^n] \sqrt{1-zA}}{[z^n]N(z)}\\*
&\sim
\Big(\frac {\gamma}{\log2 }-\frac12-\frac {\log   \pi   }{\log2 }+\frac{\log A}{2\log 2}\Big)\sqrt{A}\frac{[z^n] \sqrt{1-zA}}{-\sqrt{A}[z^n]\sqrt{1-zA}}\\
&=-\frac {\gamma}{\log2 }+\frac12+\frac {\log   \pi   }{\log2 }-\frac{\log A}{2\log 2}.
\end{align*}
The last term we consider is
\begin{align*}
\frac{2  }{\log2}&\Gamma(\chi_k)\zeta(\chi_k) A^{\frac{1-\chi_k}2}\frac{[z^n](1-zA)^{\frac{1-\chi_k}2}}{-\sqrt{A}[z^n]\sqrt{1-zA}}\\
&\sim-\frac{4 \sqrt\pi }{\log2}\frac{\Gamma(\chi_k)\zeta(\chi_k)}{\Gamma\big(\frac{\chi_k-1}{2}\big)} A^{\frac{1-\chi_k}2}n^{\chi_k/2}.
\end{align*}
Eventually we have evaluated the average value of the Horton-Strahler numbers:
\begin{theorem}
\begin{align*}
\log_4 n&- \frac{\gamma }{2\log2}- \frac{1}{\log2}+\frac32+\frac {\log   \pi   }{\log2 }-\frac{\log A}{2\log 2}
-\frac{4 \sqrt{\pi A} }{\log2}\sum_{k\neq0}\frac{\Gamma(\chi_k)\zeta(\chi_k)}{\Gamma\big(\frac{\chi_k-1}{2}\big)} A^{\frac{-\chi_k}2}n^{\chi_k/2}\\
&=\log_4 n- \frac{\gamma }{2\log2}- \frac{1}{\log2}+\frac32+\frac {\log   \pi   }{\log2 }-\frac{\log A}{2\log 2}+\psi(\log_4n),
\end{align*}
with a tiny periodic function $\psi(x)$ of period 1.
\end{theorem}	


\section{Marked ordered trees}

In \cite{Deutsch-italy} we find the following variation of ordered trees: Each rightmost edge might be marked or not, if it does not lead to an endnode (leaf).
We depict a marked edge by the red colour and draw all of them of size 4 (4 nodes):
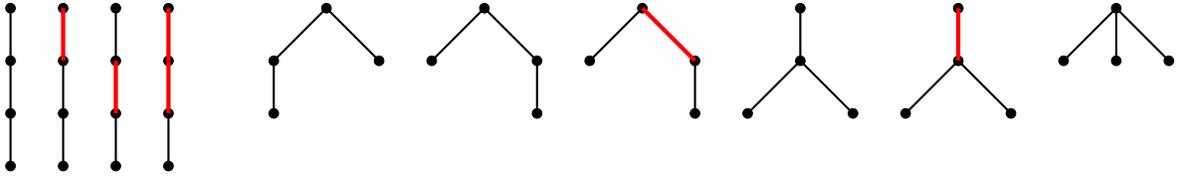
\begin{figure}[h]
	\begin{tikzpicture}[scale=0.7]

		\draw[black,fill=black] (0,0) circle (.5ex);
		\draw[black,fill=black] (0,-1) circle (.5ex);
		\draw[black,fill=black] (0,-2) circle (.5ex);
		\draw[black,fill=black] (0,-3) circle (.5ex);
		\draw [thick] (0,0) -- (0,-1)-- (0,-3) ;%
		
		\draw[black,fill=black,xshift=1cm] (0,0) circle (.5ex);
		\draw[black,fill=black,xshift=1cm] (0,-1) circle (.5ex);
		\draw[black,fill=black,xshift=1cm] (0,-2) circle (.5ex);
		\draw[black,fill=black,xshift=1cm] (0,-3) circle (.5ex);
		\draw [ultra thick,red,xshift=1cm] (0,0) -- (0,-1) ;
		\draw [thick,xshift=1cm] (0,-1) -- (0,-2) ;
		\draw [thick,xshift=1cm] (0,-2) -- (0,-3) ;

\draw[black,fill=black,xshift=2cm] (0,0) circle (.5ex);
\draw[black,fill=black,xshift=2cm] (0,-1) circle (.5ex);
\draw[black,fill=black,xshift=2cm] (0,-2) circle (.5ex);
\draw[black,fill=black,xshift=2cm] (0,-3) circle (.5ex);
\draw [thick,xshift=2cm] (0,0) -- (0,-1) ;
\draw [ultra thick,red,xshift=2cm] (0,-1) -- (0,-2) ;
\draw [thick,xshift=2cm] (0,-2) -- (0,-3) ;

\draw[black,fill=black,xshift=3cm] (0,0) circle (.5ex);
\draw[black,fill=black,xshift=3cm] (0,-1) circle (.5ex);
\draw[black,fill=black,xshift=3cm] (0,-2) circle (.5ex);
\draw[black,fill=black,xshift=3cm] (0,-3) circle (.5ex);
\draw [ultra thick,red, xshift=3cm] (0,0) -- (0,-1) ;
\draw [ultra thick,red,xshift=3cm] (0,-1) -- (0,-2) ;
\draw [thick,xshift=3cm] (0,-2) -- (0,-3) ;
		
		\draw[black,fill=black,xshift=6cm] (0,0) circle (.5ex);
		\draw[black,fill=black,xshift=6cm] (-1,-1) circle (.5ex);
		\draw[black,fill=black,xshift=6cm] (-1,-2) circle (.5ex);
		\draw[black,fill=black,xshift=6cm] (1,-1) circle (.5ex);
		\draw [thick, xshift=6cm] (0,0) -- (-1,-1) ;
		\draw [thick,xshift=6cm] (-1,-1) -- (-1,-2) ;
		\draw [thick,xshift=6cm] (0,0) -- (1,-1) ;
		
		\draw[black,fill=black,xshift=9cm] (0,0) circle (.5ex);
		\draw[black,fill=black,xshift=9cm] (-1,-1) circle (.5ex);
		\draw[black,fill=black,xshift=9cm] (1,-2) circle (.5ex);
		\draw[black,fill=black,xshift=9cm] (1,-1) circle (.5ex);
		\draw [thick, xshift=9cm] (0,0) -- (-1,-1) ;
		\draw [thick,xshift=9cm] (0,0) -- (1,-1) ;
		\draw [thick,xshift=9cm] (1,-1) -- (1,-2) ;
		
		\draw[black,fill=black,xshift=12cm] (0,0) circle (.5ex);
		\draw[black,fill=black,xshift=12cm] (-1,-1) circle (.5ex);
		\draw[black,fill=black,xshift=12cm] (1,-2) circle (.5ex);
		\draw[black,fill=black,xshift=12cm] (1,-1) circle (.5ex);
		\draw [thick, xshift=12cm] (0,0) -- (-1,-1) ;
		\draw [ultra thick,red,xshift=12cm] (0,0) -- (1,-1) ;
		\draw [thick,xshift=12cm] (1,-1) -- (1,-2) ;

\draw[black,fill=black,xshift=15cm] (0,0) circle (.5ex);
\draw[black,fill=black,xshift=15cm] (0,-1) circle (.5ex);
\draw[black,fill=black,xshift=15cm] (1,-2) circle (.5ex);
\draw[black,fill=black,xshift=15cm] (-1,-2) circle (.5ex);
\draw [thick, xshift=15cm] (0,0) -- (0,-1) ;
\draw [thick,xshift=15cm] (0,-1) -- (1,-2) ;
\draw [thick,xshift=15cm] (0,-1) -- (-1,-2) ;

\draw[black,fill=black,xshift=18cm] (0,0) circle (.5ex);
\draw[black,fill=black,xshift=18cm] (0,-1) circle (.5ex);
\draw[black,fill=black,xshift=18cm] (1,-2) circle (.5ex);
\draw[black,fill=black,xshift=18cm] (-1,-2) circle (.5ex);
\draw [ultra thick, red, xshift=18cm] (0,0) -- (0,-1) ;
\draw [ thick,xshift=18cm] (0,-1) -- (1,-2) ;
\draw [thick,xshift=18cm] (0,-1) -- (-1,-2) ;
		
		\draw[black,fill=black,xshift=21cm] (0,0) circle (.5ex);
		\draw[black,fill=black,xshift=21cm] (0,-1) circle (.5ex);
		\draw[black,fill=black,xshift=21cm] (1,-1) circle (.5ex);
		\draw[black,fill=black,xshift=21cm] (-1,-1) circle (.5ex);
		\draw [thick, xshift=21cm] (0,0) -- (0,-1) ;
		\draw [thick,xshift=21cm] (0,0) -- (1,-1) ;
		\draw [thick,xshift=21cm] (0,0) -- (-1,-1) ;
		
		\end{tikzpicture}
	\caption{All 10 marked ordered trees with 4 nodes.}
\end{figure}

Now we move to a symbolic equation for the marked ordered trees:
\begin{figure}[h]\small
	\begin{tikzpicture}[scale=1.0,
		s1/.style={circle=10pt,draw=black!90,thick},
		s2/.style={rectangle,draw=black!50,thick},scale=0.5]
		
		\node at ( -3.0,0) { $\mathscr{A}$};
		
		\node at (-1.5,0) { $=$};
		\node(c) at (-0.4,0)[s1]{};
		\node at (1.3,0) {$+$};

		\node(d) at (5,1)[s1]{};
		\node(e) at (3,-1){ $\mathscr{A}$};
\node(ee) at (5,-1){$\cdots$};
		\node(f) at (7,-1){ $\mathscr{A}$};
		\path [draw,-,black!90] (d) -- (e) node{};%
		\path [draw,-,black!90] (d) -- (f) node{};%
		\node(g) at (9,-1)[s1]{ };
		\path [draw,-,black,ultra thick] (d) -- (g) node {};%
		
		\node[xshift=5cm] at (0.7,0) {$+$};

		\node[xshift=5cm](d) at (5,1)[s1]{};
		\node[xshift=5cm](e) at (3,-1){ $\mathscr{A}$};
		\node[xshift=5cm](ee) at (5,-1){$\cdots$};
		\node[xshift=5cm](f) at (7,-1){ $\mathscr{A}$};
		\path [draw,-,black!90] (d) -- (e) node{};%
		\path [draw,-,black!90] (d) -- (f) node{};%
		\node[xshift=5cm](g) at (9,-1){ $\mathscr{A}\setminus\{\circ\}$};
		\path [draw,-,black,ultra thick] (d) -- (8.8+10,-0.5) node {};%
		
		\node[xshift=10cm] at (0.7,0) {$+$};

		\node[xshift=10cm](d) at (5,1)[s1]{};
		\node[xshift=10cm](e) at (3,-1){ $\mathscr{A}$};
		\node[xshift=10cm](ee) at (5,-1){$\cdots$};
		\node[xshift=10cm](f) at (7,-1){ $\mathscr{A}$};
		\path [draw,-,black!90] (d) -- (e) node{};%
		\path [draw,-,black!90] (d) -- (f) node{};%
		\node[xshift=10cm](g) at (9,-1){ $\mathscr{A}\setminus\{\circ\}$};
		\path [draw,-,black,red,ultra thick] (d) -- (8.8+20,-0.5) node {};%

	\end{tikzpicture}
	\caption{Symbolic equation for marked ordered trees.\\ $\mathscr{A}\cdots\mathscr{A}$ refers to $\ge0$ copies of $\mathscr{A}$.}
\end{figure}

In terms of generating functions,
\begin{equation*}
A=z+\frac{z}{1-A}z+\frac{z}{1-A}2(A-z),
\end{equation*}
with the solution
\begin{equation*}
A(z)=\frac{1-z-\sqrt{1-6z+5z^2}}{2}=z+z^2+z^3+3z^3+10z^4+36z^5+\cdots.
\end{equation*}

The importance of this family of trees lies in the bijection to skew Dyck paths, as given in \cite{Deutsch-italy}. One walks around the tree as one usually does and translates it into a Dyck path.
The only difference are the red edges. On the way down, nothing special is to be reported, but on the way up, it is translated into  a skew step $(-1,-1)$. The present author believes that trees are more manageable when it comes to enumeration issues than skew Dyck paths.

The 10 trees of Figure 1 translate as follows:
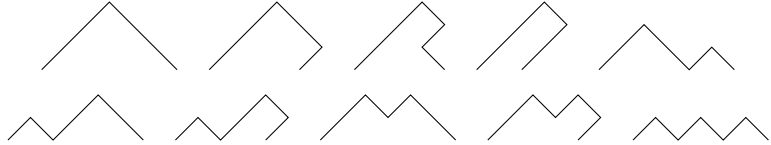
\begin{equation*}
		\begin{tikzpicture}[scale=0.3]
		\draw (0,0)--(3,3)--(6,0);
	\end{tikzpicture}
\quad
\begin{tikzpicture}[scale=0.3]
	\draw (0,0)--(3,3)--(5,1)--(4,0);
\end{tikzpicture}
\quad
\begin{tikzpicture}[scale=0.3]
	\draw (0,0)--(3,3)--(4,2)--(3,1)--(4,0);
\end{tikzpicture}
\quad
\begin{tikzpicture}[scale=0.3]
	\draw (0,0)--(3,3)--(4,2)--(3,1)--(2,0);
\end{tikzpicture}
\quad
\begin{tikzpicture}[scale=0.3]
	\draw (0,0)--(2,2)--(4,0)--(5,1)--(6,0);
\end{tikzpicture}
\end{equation*}
\begin{equation*}
\begin{tikzpicture}[scale=0.3]
	\draw (0,0)--(1,1)--(2,0)--(4,2)--(6,0);
\end{tikzpicture}
\quad
\begin{tikzpicture}[scale=0.3]
	\draw (0,0)--(1,1)--(2,0)--(4,2)--(5,1)--(4,0);
\end{tikzpicture}
\quad
\begin{tikzpicture}[scale=0.3]
	\draw (0,0)--(1,1)--(2,2)--(3,1)--(4,2)--(6,0);
\end{tikzpicture}
\quad
\begin{tikzpicture}[scale=0.3]
	\draw (0,0)--(1,1)--(2,2)--(3,1)--(4,2)--(5,1)--(4,0);
\end{tikzpicture}
\quad
\begin{tikzpicture}[scale=0.3]
	\draw (0,0)--(1,1)--(2,0)--(3,1)--(4,0)--(5,1)--(6,0);
\end{tikzpicture}
\end{equation*}

\section{Parameters of marked ordered trees}

There are many parameters, usually considered in the context of ordered trees, that can be considered 
for marked ordered trees. Of course, we cannot be encyclopedic, also, to keep the balance between the
other structures that are considered in this paper. We just consider a few parameters and leave further
analysis to the future.

\subsection*{The number of leaves}

To get this, it is most natural to use an additional variable $u$ when translating the symbolic equation, so
that $z^nu^k$ refers to trees with $n$ nodes and $k$ leaves. One obtains
\begin{equation*}
F=zu+\frac{z}{1-F}\bigl(zu2(F-zu)\bigr),
\end{equation*}
with the solution
\begin{align*}
F(z,u)&=-z+\frac{zu}2+\frac12-\frac12\sqrt {4{z}^{2}-4z+{z}^{2}{u}^{2}-2zu+1}\\*
&=zu+{z}^{2}u+ \left( 2u+{u}^{2}\right) {z}^{3}+ \left( 4u+5{u}^{2}+{u}^{3} \right) {z}^{4}+\cdots.
\end{align*}
The factor $4u+5{u}^{2}+{u}^{3}$ corresponds to the 10 trees in Figure 1.

Of interest is also the average number of leaves, when all marked ordered trees of size $n$ are considered to be 
equally likely. For that, we differentiate $F(z,u)$ w. r. t. $u$, followed by $u:=1$, with the result
\begin{equation*}
\frac{z}{2}+\frac{z-z^2}{2\sqrt{1-6z+5z^2}}=\frac{z}{1-v}, \quad\text{with the usual}\quad z=\frac{v}{1+3v+v^2}.
\end{equation*}
Since $F(z,1)=z(1+v)$, it follows that the average is asymptotic to 
\begin{align*}
\frac{[z^{n+1}]\frac{z}{1-v}}{[z^{n+1}]z(1+v)}&=\frac{[z^{n}]\frac{1}{1-v}}{[z^{n}](1+v)}=\frac{[z^n]\frac1{\sqrt5}\frac1{\sqrt{1-5z}}}{5^{n+\frac12}\frac1{2\sqrt\pi}n^{3/2}}\\
&=\frac{\frac{n^{-1/2}}{\Gamma(\frac12)}}{5^{n+\frac12}\frac1{2\sqrt\pi}n^{3/2}}=\frac{n}{10}.
\end{align*}
Note that the corresponding number for ordered trees (unmarked) is $\frac n2$, so we have significantly less leaves here.

\subsection*{The height}

As in the seminal paper \cite{BrKnRi72}, we define the height in terms of the longest chain of nodes from the root to a leaf. 
Further, let $p_h=p_h(z)$ be the generating function of marked ordered trees of height $\le h$. From the symbolic equation,
\begin{align*}
p_{h+1}=z+\frac{z^2}{1-p_h}+\frac{2z^2(p_h-z)}{1-p_h}=-z+\frac{2z-z^2}{1-p_h},\quad h\ge1,\ p_1=z.
\end{align*}
By some creative guessing, separating numerator and denominator, we find the solution
\begin{equation*}
p_h=z(1+v)\frac{(1+2v)^{h-1}-v^h(v+2)^{h-1}}{(1+2v)^{h-1}-v^{h+1}(v+2)^{h-1}},
\end{equation*}
which is easy to prove by induction. The limit for $h\to\infty$ is $z(1+v)$, the generating function of \emph{all} marked ordered trees, as expected.
Taking differences, we get the generating functions of trees of height $>h$:
\begin{align*}
z(1+v)&-z(1+v)\frac{(1+2v)^{h-1}-v^h(v+2)^{h-1}}{(1+2v)^{h-1}-v^{h+1}(v+2)^{h-1}}\\
&=z(1+v)\frac{(1+2v)^{h-1}-v^{h+1}(v+2)^{h-1}-(1+2v)^{h-1}+v^h(v+2)^{h-1}}{(1+2v)^{h-1}-v^{h+1}(v+2)^{h-1}}\\
&=z(1-v^2)\frac{(v+2)^{h-1}v^h}{(1+2v)^{h-1}-v^{h+1}(v+2)^{h-1}}.
\end{align*}
From this, the average height can be worked out, as in the model paper \cite{HPW}. We sketch the essential steps.
For the average height, one needs
\begin{equation*}
\sum_{h\ge0}z(1-v^2)\frac{(v+2)^{h-1}v^h}{(1+2v)^{h-1}-v^{h+1}(v+2)^{h-1}}
\end{equation*}
and its behaviour around $v=1$, viz.
\begin{equation*}
2z(1-v)\sum_{h\ge0}\frac{3^{h-1}v^h}{3^{h-1}-v^{h+1}3^{h-1}}
\sim2z(1-v)\sum_{h\ge1}\frac{v^h}{1-v^{h}}.
\end{equation*}
The behaviour of the series can be taken straight from \cite{HPW}.

We find there
\begin{equation*}
\sum_{h\ge1}\frac{v^h}{1-v^{h}}=-\frac{\log(1-v)}{1-v},
\end{equation*}
and further
\begin{equation*}
	\sum_{h\ge0}z(1-v^2)\frac{(v+2)^{h-1}v^h}{(1+2v)^{h-1}-v^{h+1}(v+2)^{h-1}}
\sim-	2z\log(1-v),
\end{equation*}
so that the coefficient of $z^{n+1}$ is asymptotic to $-2[z^n]\log(1-v)$. Since $1-v\sim \sqrt5\sqrt{1-5z}$,
\begin{equation*}
-	2z\log(1-v)\sim -2z\log\sqrt{1-5z}= -z\log(1-5z),
\end{equation*}
and the coefficient of $z^{n+1}$ in it is asymptotic to $\frac{5^n}{n}$. This has to be divided (as derived earlier) by
\begin{equation*}
5^{n+\frac12}\frac{1}{2\sqrt\pi n^{3/2}},
\end{equation*}
with the result
\begin{equation*}
2\frac{5^n}{n}\frac1{5^{n+\frac12}}\sqrt\pi n^{3/2}=\frac{2}{\sqrt5}\sqrt{\pi n}.
\end{equation*}
Note that the constant in front of $\sqrt{\pi n}$ for ordered trees is $\frac{2}{\sqrt4}=1$, so the average height for marked ordered trees is indeed a bit smaller
thanks to the extra markings.

\section{A bijection between multi-edge trees and 3-coloured Motzkin paths}

Multi-edge trees are like ordered (planar, plane, \dots) trees, but instead of edges there are multiple edges. When drawing such a tree,
instead of drawing, say 5 parallel edges, we just draw one edge and put the number 5 on it as a label. These trees were studied in
\cite{polish, HPW}. For the enumeration, one must count edges. The generating function $F(z)$ satisfies
\begin{equation*}
	F(z)=\sum_{k\ge0}\Big(\frac{z}{1-z}F(z)\Big)^k=\frac1{1-\frac{z}{1-z}F(z)},
\end{equation*}
whence
\begin{equation*}
	F(z)=\frac{1-z-\sqrt{1-6z+5z^2}}{2z}=1+z+3{z}^{2}+10{z}^{3}+36{z}^{4}+137{z}^{5}+543{z}^{6}+\cdots.
\end{equation*}
The coefficients form once again sequence A002212 in \cite{OEIS}.

A Motzkin path consists of up-steps, down-steps, and horizontal steps, see sequence A091965 in \cite{OEIS} and the references given there. As Dyck paths, they start at the origin and end, after $n$ steps again at the $x$-axis, but are not allowed to go below the $x$-axis. A 3-coloured Motzkin path is built as a Motzkin path, but there
are 3 different types of horizontal steps, which we call \emph{red, green, blue}. The generating function $M(z)$ satisfies
\begin{equation*}
	M(z)=1+3zM(z)+z^2M(z)^2=\frac{1-3z-\sqrt{1-6z+5z^2}}{2z^2}, \quad\text{or}\quad F(z)=1+zM(z).
\end{equation*}
So multi-edge trees with $N$ edges (counting the multiplicities) correspond to  3-coloured Motzkin paths of length $N-1$.

The purpose of this note is to describe a bijection. It transforms trees into paths, but all steps are reversible.

\subsection*{The details}

As a first step, the multiplicities will be ignored, and the tree then has only $n$ edges. The standard translation of such tree into the world of Dyck paths,
which is in every book on combinatorics, leads to a Dyck path of length $2n$.
Then the Dyck path will transformed bijectively to a 2-coloured Motzkin path of length $n-1$ (the colours used are red and green).
This transformation plays a prominent role in \cite{Shapiro}, but is most likely much older. I believe that people like Viennot know this for 40 years.
I would be glad to get a proper historic account from the gentle readers.

The last step is then to use the third colour (blue) to deal with the multiplicities.

The first up-step and the last down-step of the Dyck path will be deleted. Then, the remaining $2n-2$ steps are coded pairwise into a 2-Motzkin path of length $n-1$:
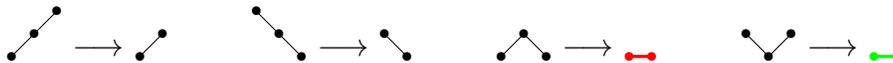
\begin{equation*}
	\begin{tikzpicture}[scale=0.3]
		\path (0,0) node(x1) {\tiny$\bullet$} ;
		\path (1,1) node(x2) {\tiny$\bullet$};
		\path (2,2) node(x3) {\tiny$\bullet$};
		
		\draw (0,0) -- (2,2);

	\end{tikzpicture}
	\raisebox{0.5 em}{$\longrightarrow$}
	\begin{tikzpicture}[scale=0.3]
		\path (0,0) node(x1) {\tiny$\bullet$} ;
		\path (1,1) node(x2) {\tiny$\bullet$};
		
		\draw (0,0) -- (1,1);

	\end{tikzpicture}
	\qquad
	\begin{tikzpicture}[scale=0.3]
		\path (0,2) node(x1) {\tiny$\bullet$} ;
		\path (1,1) node(x2) {\tiny$\bullet$};
		\path (2,0) node(x3) {\tiny$\bullet$};
		
		\draw (0,2) -- (2,0);

	\end{tikzpicture}
	\raisebox{0.5 em}{$\longrightarrow$}
	\begin{tikzpicture}[scale=0.3]
		\path (0,1) node(x1) {\tiny$\bullet$} ;
		\path (1,0) node(x2) {\tiny$\bullet$};
		
		\draw (0,1) -- (1,0);

	\end{tikzpicture}
	\qquad
	\begin{tikzpicture}[scale=0.3]
		\path (0,0) node(x1) {\tiny$\bullet$} ;
		\path (1,1) node(x2) {\tiny$\bullet$};
		\path (2,0) node(x3) {\tiny$\bullet$};
		
		\draw (0,0) -- (1,1) -- (2,0);

	\end{tikzpicture}
	\raisebox{0.5 em}{$\longrightarrow$}
	\begin{tikzpicture}[scale=0.3]
		\path (0,0) node[red](x1) {\tiny$\bullet$} ;
		\path (1,0) node[red](x2) {\tiny$\bullet$};
		
		\draw[red, very thick] (0,0) -- (1,0);

	\end{tikzpicture}
	\qquad
	\begin{tikzpicture}[scale=0.3]
		\path (0,0) node(x1) {\tiny$\bullet$} ;
		\path (1,-1) node(x2) {\tiny$\bullet$};
		\path (2,0) node(x3) {\tiny$\bullet$};
		
		\draw (0,0) -- (1,-1) -- (2,0);

	\end{tikzpicture}
	\raisebox{0.5 em}{$\longrightarrow$}
	\begin{tikzpicture}[scale=0.3]
		\path (0,0) node[green](x1) {\tiny$\bullet$} ;
		\path (1,0) node[green](x2) {\tiny$\bullet$};
		
		\draw[green, very thick] (0,0) -- (1,0);

	\end{tikzpicture}
\end{equation*}

The last step is to deal with the multiplicities. If an edge is labelled with the number $a$, we will insert $a-1$ horizontal blue steps in the following way:
Since there are currently $n-1$ symbols in the path, we have $n$ possible positions to enter something (in the beginning, in the end, between symbols).
We go through the tree in pre-order, and enter the multiplicities one by one using the blue horizontal steps.

To make this procedure more clear, we prepared a list of 10 multi-edge trees with 3 edges, and the corresponding 3-Motzkin paths of length 2, with intermediate steps 
completely worked out:

	\begin{center}
		\begin{table}[h]
			\begin{tabular}{c | c | c  |c}
				\text{Multi-edge tree  }&\text{Dyck path}&\text{2-Motzkin path}&\text{blue edges added}\\
				\hline\hline
				\begin{tikzpicture}[scale=0.5]
					\path (0,0) node(x1) {\tiny$\bullet$} ;
					\path (0,-1) node(x2) {\tiny$\bullet$};
					\path (0,-2) node(x3) {\tiny$\bullet$};
					\path (0,-3) node(x4) {\tiny$\bullet$};
					\draw (0,0) -- (0,-1)node[pos=0.5,left]{\tiny1} ;
					\draw (0,-1) -- (0,-2)node[pos=0.5,left]{\tiny1} ;
					\draw (0,-2) -- (0,-3)node[pos=0.5,left]{\tiny 1} ;

				\end{tikzpicture}
				& \begin{tikzpicture}[scale=0.45]
					
					\draw (0,0) -- (3,3) --(6,0);

				\end{tikzpicture}
				& 
				\begin{tikzpicture}[scale=0.45]
					
					\draw[thick] (0,0) -- (1,1) --(2,0);

				\end{tikzpicture}
				& \begin{tikzpicture}[scale=0.45]
					
					\draw[thick] (0,0) -- (1,1) --(2,0);

				\end{tikzpicture}\\
				
				\hline

				\begin{tikzpicture}[scale=0.5]
					\path (0,0) node(x1) {\tiny$\bullet$} ;
					\path (0,-1) node(x2) {\tiny$\bullet$};
					\path (0,-2) node(x3) {\tiny$\bullet$};
					
					\draw (0,0) -- (0,-1)node[pos=0.5,left]{\tiny2} ;
					\draw (0,-1) -- (0,-2)node[pos=0.5,left]{\tiny1} ;

				\end{tikzpicture}
				& \begin{tikzpicture}[scale=0.45]
					
					\draw (0,0) -- (2,2) --(4,0);

				\end{tikzpicture}
				& \begin{tikzpicture}[scale=0.45]
					
					\draw [red,thick](0,0) -- (1,0);

				\end{tikzpicture}
				& \begin{tikzpicture}[scale=0.45]
					\draw [blue,thick](0,0) -- (1,0);
					\draw [red,thick](1,0) -- (2,0);

				\end{tikzpicture}\\
				
				\hline
				\begin{tikzpicture}[scale=0.5]
					\path (0,0) node(x1) {\tiny$\bullet$} ;
					\path (0,-1) node(x2) {\tiny$\bullet$};
					\path (0,-2) node(x3) {\tiny$\bullet$};
					
					\draw (0,0) -- (0,-1)node[pos=0.5,left]{\tiny1} ;
					\draw (0,-1) -- (0,-2)node[pos=0.5,left]{\tiny2} ;

				\end{tikzpicture}
				& \begin{tikzpicture}[scale=0.45]
					
					\draw (0,0) -- (2,2) --(4,0);

				\end{tikzpicture}& \begin{tikzpicture}[scale=0.45]
					
					\draw [red,thick](0,0) -- (1,0);

				\end{tikzpicture}& \begin{tikzpicture}[scale=0.45]
					\draw [red,thick](0,0) -- (1,0);
					\draw [blue,thick](1,0) -- (2,0);

				\end{tikzpicture}\\
				
				\hline
				\begin{tikzpicture}[scale=0.5]
					\path (0,0) node(x1) {\tiny$\bullet$} ;
					\path (0,-1) node(x2) {\tiny$\bullet$};

					\draw (0,0) -- (0,-1)node[pos=0.5,left]{\tiny3} ;

				\end{tikzpicture}
				& \begin{tikzpicture}[scale=0.45]
					
					\draw (0,0) -- (1,1) --(2,0);

				\end{tikzpicture} & & \begin{tikzpicture}[scale=0.45]
					\draw [blue,thick](0,0) -- (2,0);

				\end{tikzpicture}\\
				
				\hline 
				\begin{tikzpicture}[scale=0.5]
					\path (0,0) node(x1) {\tiny$\bullet$} ;
					\path (-1,-1) node(x2) {\tiny$\bullet$};
					\path (-1,-2) node(x3) {\tiny$\bullet$};
					\path (1,-1) node(x4) {\tiny$\bullet$};
					\draw (0,0) -- (-1,-1)node[pos=0.3,left]{\tiny1} ;
					\draw (-1,-1) -- (-1,-2)node[pos=0.3,left]{\tiny1} ;
					\draw (0,0) -- (1,-1)node[pos=0.3,right]{\tiny1} ;

				\end{tikzpicture}
				& \begin{tikzpicture}[scale=0.45]
					
					\draw (0,0) -- (2,2) --(4,0)--(5,1)--(6,0);

				\end{tikzpicture}& \begin{tikzpicture}[scale=0.45]
					
					\draw [red,thick](0,0) -- (1,0);
					\draw [green,thick](1,0) -- (2,0);
					
				\end{tikzpicture}& \begin{tikzpicture}[scale=0.45]
					
					\draw [red,thick](0,0) -- (1,0);
					\draw [green,thick](1,0) -- (2,0);
					
				\end{tikzpicture}\\
				
				\hline
				\begin{tikzpicture}[scale=0.5]
					\path (0,0) node(x1) {\tiny$\bullet$} ;
					\path (-1,-1) node(x2) {\tiny$\bullet$};
					
					\path (1,-1) node(x4) {\tiny$\bullet$};
					\draw (0,0) -- (-1,-1)node[pos=0.3,left]{\tiny2} ;
					
					\draw (0,0) -- (1,-1)node[pos=0.3,right]{\tiny1} ;

				\end{tikzpicture}
				& \begin{tikzpicture}[scale=0.45]
					
					\draw (0,0) -- (1,1) --(2,0)--(3,1)--(4,0);

				\end{tikzpicture}& \begin{tikzpicture}[scale=0.45]

					\draw [green,thick](0,0) -- (1,0);
					
				\end{tikzpicture}			& \begin{tikzpicture}[scale=0.45]
					
					\draw [blue,thick](0,0) -- (1,0);
					\draw [green,thick](1,0) -- (2,0);
					
				\end{tikzpicture}\\
				
				\hline
				\begin{tikzpicture}[scale=0.5]
					\path (0,0) node(x1) {\tiny$\bullet$} ;
					\path (-1,-1) node(x2) {\tiny$\bullet$};
					
					\path (1,-1) node(x4) {\tiny$\bullet$};
					\draw (0,0) -- (-1,-1)node[pos=0.3,left]{\tiny1} ;
					
					\draw (0,0) -- (1,-1)node[pos=0.3,right]{\tiny2} ;

				\end{tikzpicture}
				& \begin{tikzpicture}[scale=0.45]
					
					\draw (0,0) -- (1,1) --(2,0)--(3,1)--(4,0);

				\end{tikzpicture}& \begin{tikzpicture}[scale=0.45]

					\draw [green,thick](0,0) -- (1,0);
					
				\end{tikzpicture}& \begin{tikzpicture}[scale=0.45]
					
					\draw [green,thick](0,0) -- (1,0);
					\draw [blue,thick](1,0) -- (2,0);
					
				\end{tikzpicture}\\
				
				\hline 
				\begin{tikzpicture}[scale=0.5]
					\path (-1,0) node(x1) {\tiny$\bullet $} ;
					\path (-1,1) node(x2) {\tiny$\bullet$};
					\path (-2,2) node(x3) {\tiny$\bullet$};
					\path (-3,1) node(x4) {\tiny$\bullet$};
					\draw (-1,0) -- (-1,1)node[pos=0.7,right]{\tiny1} ;
					\draw (-1,1) -- (-2,2)node[pos=0.7,right]{\tiny1} ;
					\draw (-2,2) -- (-3,1)node[pos=0.3,left]{\tiny1} ;

				\end{tikzpicture}
				& \begin{tikzpicture}[scale=0.45]
					
					\draw (0,0) -- (1,1) --(2,0)--(4,2)--(6,0);

				\end{tikzpicture}& \begin{tikzpicture}[scale=0.45]

					\draw [green,thick](0,0) -- (1,0);
					\draw[red,thick](1,0)--(2,0);			
				\end{tikzpicture}& \begin{tikzpicture}[scale=0.45]

					\draw [green,thick](0,0) -- (1,0);
					\draw[red,thick](1,0)--(2,0);			
				\end{tikzpicture}\\
				
				\hline
				\begin{tikzpicture}[scale=0.5]
					\path (0,0) node(x1) {\tiny$\bullet$} ;
					\path (-1,-1) node(x2) {\tiny$\bullet$};
					\path (0,-1) node(x3) {\tiny$\bullet$};
					\path (1,-1) node(x4) {\tiny$\bullet$};
					\draw (0,0) -- (-1,-1)node[pos=0.3,left]{\tiny1} ;
					\draw (0,0) -- (0,-1)node[pos=0.6]{\tiny\;\;1} ;
					\draw (0,0) -- (1,-1)node[pos=0.3,right]{\tiny1} ;

				\end{tikzpicture}
				& \begin{tikzpicture}[scale=0.45]
					
					\draw (0,0) -- (1,1) --(2,0)--(3,1)--(4,0)--(5,1)--(6,0);

				\end{tikzpicture}& \begin{tikzpicture}[scale=0.45]

					\draw [green,thick](0,0) -- (1,0);
					\draw [green,thick](1,0) -- (2,0);
					
				\end{tikzpicture}& \begin{tikzpicture}[scale=0.45]

					\draw [green,thick](0,0) -- (1,0);
					\draw [green,thick](1,0) -- (2,0);
					
				\end{tikzpicture}\\
				
				\hline 
				\begin{tikzpicture}[scale=0.5]
					\path (0,0) node(x1) {\tiny$\bullet$} ;
					\path (0,-1) node(x2) {\tiny$\bullet$};
					\path (-1,-2) node(x3) {\tiny$\bullet$};
					\path (1,-2) node(x4) {\tiny$\bullet$};
					\draw (0,0) -- (0,-1)node[pos=0.5,left]{\tiny1} ;
					\draw (0,-1) -- (1,-2)node[pos=0.5,right]{\tiny1} ;
					\draw (0,-1) -- (-1,-2)node[pos=0.5,left]{\tiny 1} ;

				\end{tikzpicture}
				& \begin{tikzpicture}[scale=0.45]
					
					\draw (0,0) -- (2,2) --(3,1)--(4,2)--(6,0);

				\end{tikzpicture}& \begin{tikzpicture}[scale=0.45]

					\draw [red,thick](0,0) -- (1,0);
					\draw [red,thick](1,0) -- (2,0);
					
				\end{tikzpicture}& \begin{tikzpicture}[scale=0.45]

					\draw [red,thick](0,0) -- (1,0);
					\draw [red,thick](1,0) -- (2,0);
					
				\end{tikzpicture}\\

			\end{tabular}
			
			\caption{First row is a multi-edge tree with 3 edges, second row is the standard Dyck path (multiplicities ignored), third row is cutting off first and last step, and then translated pairs of steps, fourth row is inserting blued horizontal edges, according to multiplicities.}
		\end{table}	
	\end{center}

	\subsection*{Connecting unary-binary trees with multi-edge trees} 
	
	This is not too difficult: We start from multi-edge trees, and ignore the multiplicities at the moment. Then we apply the classical rotation correspondence (also called: natural correspondence).
	Then we add vertical edges, if the multiplicity is higher than 1. To be precise, if there is a node, and an edge with multiplicity $a$ leads to it from the top, we insert $a-1$ extra nodes in a chain on the top, and connect them with unary branches. The following example with 10 objects will help to understand this procedure.
	After that, all the structures studied in this paper are connected with bijections.
	
\begin{center}
	\begin{table}[h]
		\begin{tabular}{c | c | c  |c}
			\text{Multi-edge tree  }&\text{Binary tree (rotation)}&\text{vertical edges added}\\
			\hline\hline
			\begin{tikzpicture}[scale=0.5]
				\path (0,0) node(x1) {\tiny$\bullet$} ;
				\path (0,-1) node(x2) {\tiny$\bullet$};
				\path (0,-2) node(x3) {\tiny$\bullet$};
				\path (0,-3) node(x4) {\tiny$\bullet$};
				\draw (0,0) -- (0,-1)node[pos=0.5,left]{\tiny1} ;
				\draw (0,-1) -- (0,-2)node[pos=0.5,left]{\tiny1} ;
				\draw (0,-2) -- (0,-3)node[pos=0.5,left]{\tiny 1} ;

			\end{tikzpicture}

			& 			\begin{tikzpicture}[scale=0.5]
				\path (0,0) node(x1) {\tiny$\bullet$} ;
				\path (-1,-1) node(x2) {\tiny$\bullet$};
				\path (-2,-2) node(x3) {\tiny$\bullet$};
				
				\draw (0,0) -- (-1,-1);
				\draw (-1,-1) -- (-2,-2) ;

			\end{tikzpicture}
			
			& \begin{tikzpicture}[scale=0.5]
				\path (0,0) node(x1) {\tiny$\bullet$} ;
				\path (-1,-1) node(x2) {\tiny$\bullet$};
				\path (-2,-2) node(x3) {\tiny$\bullet$};
				
				\draw (0,0) -- (-1,-1);
				\draw (-1,-1) -- (-2,-2) ;

			\end{tikzpicture}
			\\
			\hline 
			\begin{tikzpicture}[scale=0.5]
				\path (0,0) node(x1) {\tiny$\bullet$} ;
				\path (0,-1) node(x2) {\tiny$\bullet$};
				\path (0,-2) node(x3) {\tiny$\bullet$};
				
				\draw (0,0) -- (0,-1)node[pos=0.5,left]{\tiny2} ;
				\draw (0,-1) -- (0,-2)node[pos=0.5,left]{\tiny1} ;

			\end{tikzpicture}
			& \begin{tikzpicture}[scale=0.5]
				\path (0,0) node(x1) {\tiny$\bullet$} ;
				\path (-1,-1) node(x2) {\tiny$\bullet$};
				
				\draw (0,0) -- (-1,-1);

			\end{tikzpicture}
			&\begin{tikzpicture}[scale=0.5]
				\path (0,0) node(x1) {\tiny$\bullet$} ;
				\path (-1,-1) node(x2) {\tiny$\bullet$};
				\path (0,1) node(x3) {\tiny$\bullet$};
				
				\draw (0,0) -- (-1,-1);
				\draw (0,0) -- (0,1) ;

			\end{tikzpicture}
			\\
			\hline 
			\begin{tikzpicture}[scale=0.5]
				\path (0,0) node(x1) {\tiny$\bullet$} ;
				\path (0,-1) node(x2) {\tiny$\bullet$};
				\path (0,-2) node(x3) {\tiny$\bullet$};
				
				\draw (0,0) -- (0,-1)node[pos=0.5,left]{\tiny1} ;
				\draw (0,-1) -- (0,-2)node[pos=0.5,left]{\tiny2} ;

			\end{tikzpicture}
			& \begin{tikzpicture}[scale=0.5]
				\path (0,0) node(x1) {\tiny$\bullet$} ;
				\path (-1,-1) node(x2) {\tiny$\bullet$};
				
				\draw (0,0) -- (-1,-1);

			\end{tikzpicture}
			&\begin{tikzpicture}[scale=0.5]
				\path (-1,0) node(x1) {\tiny$\bullet$} ;
				\path (-1,-1) node(x2) {\tiny$\bullet$};
				\path (0,1) node(x3) {\tiny$\bullet$};
				
				\draw (-1,0) -- (-1,-1);
				\draw (-1,0) -- (0,1) ;

			\end{tikzpicture}
			\\
			\hline
			\begin{tikzpicture}[scale=0.5]
				\path (0,0) node(x1) {\tiny$\bullet$} ;
				\path (0,-1) node(x2) {\tiny$\bullet$};

				\draw (0,0) -- (0,-1)node[pos=0.5,left]{\tiny3} ;

			\end{tikzpicture}
			& \begin{tikzpicture}[scale=0.5]
				\path (0,0) node(x1) {\tiny$\bullet$} ;

			\end{tikzpicture}
			&\begin{tikzpicture}[scale=0.5]
				\path (0,0) node(x1) {\tiny$\bullet$} ;
				\path (0,-1) node(x2) {\tiny$\bullet$};
				\path (0,1) node(x3) {\tiny$\bullet$};
				
				\draw (0,0) -- (0,-1);
				\draw (0,0) -- (0,1) ;

			\end{tikzpicture}
			\\
			
			\hline 
			\begin{tikzpicture}[scale=0.5]
				\path (0,0) node(x1) {\tiny$\bullet$} ;
				\path (-1,-1) node(x2) {\tiny$\bullet$};
				\path (-1,-2) node(x3) {\tiny$\bullet$};
				\path (1,-1) node(x4) {\tiny$\bullet$};
				\draw (0,0) -- (-1,-1)node[pos=0.3,left]{\tiny1} ;
				\draw (-1,-1) -- (-1,-2)node[pos=0.3,left]{\tiny1} ;
				\draw (0,0) -- (1,-1)node[pos=0.3,right]{\tiny1} ;

			\end{tikzpicture}
			& \begin{tikzpicture}[scale=0.5]
				\path (0,-2) node(x1) {\tiny$\bullet$} ;
				\path (-1,-1) node(x2) {\tiny$\bullet$};
				\path (-2,-2) node(x3) {\tiny$\bullet$};
				
				\draw (0,-2) -- (-1,-1);
				\draw (-1,-1) -- (-2,-2) ;

			\end{tikzpicture}
			& \begin{tikzpicture}[scale=0.5]
				\path (0,-2) node(x1) {\tiny$\bullet$} ;
				\path (-1,-1) node(x2) {\tiny$\bullet$};
				\path (-2,-2) node(x3) {\tiny$\bullet$};
				
				\draw (0,-2) -- (-1,-1);
				\draw (-1,-1) -- (-2,-2) ;

			\end{tikzpicture}
			\\
			\hline
			\begin{tikzpicture}[scale=0.5]
				\path (0,0) node(x1) {\tiny$\bullet$} ;
				\path (-1,-1) node(x2) {\tiny$\bullet$};
				
				\path (1,-1) node(x4) {\tiny$\bullet$};
				\draw (0,0) -- (-1,-1)node[pos=0.3,left]{\tiny2} ;
				
				\draw (0,0) -- (1,-1)node[pos=0.3,right]{\tiny1} ;

			\end{tikzpicture}
			& \begin{tikzpicture}[scale=0.5]
				\path (0,-2) node(x1) {\tiny$\bullet$} ;
				\path (-1,-1) node(x2) {\tiny$\bullet$};
				
				\draw (0,-2) -- (-1,-1);

			\end{tikzpicture}
			&  	\begin{tikzpicture}[scale=0.5]
				\path (0,0) node(x1) {\tiny$\bullet$} ;
				\path (1,-1) node(x2) {\tiny$\bullet$};
				\path (0,1) node(x3) {\tiny$\bullet$};
				
				\draw (0,0) -- (1,-1);
				\draw (0,0) -- (0,1) ;

			\end{tikzpicture}
			\\
			\hline
			\begin{tikzpicture}[scale=0.5]
				\path (0,0) node(x1) {\tiny$\bullet$} ;
				\path (-1,-1) node(x2) {\tiny$\bullet$};
				
				\path (1,-1) node(x4) {\tiny$\bullet$};
				\draw (0,0) -- (-1,-1)node[pos=0.3,left]{\tiny1} ;
				
				\draw (0,0) -- (1,-1)node[pos=0.3,right]{\tiny2} ;

			\end{tikzpicture}
			& \begin{tikzpicture}[scale=0.5]
				\path (0,-2) node(x1) {\tiny$\bullet$} ;
				\path (-1,-1) node(x2) {\tiny$\bullet$};
				
				\draw (0,-2) -- (-1,-1);

			\end{tikzpicture}
			
			&  \begin{tikzpicture}[scale=0.5]
				\path (0,0) node(x1) {\tiny$\bullet$} ;
				\path (0,-1) node(x2) {\tiny$\bullet$};
				\path (-1,1) node(x3) {\tiny$\bullet$};
				
				\draw (0,0) -- (-1,1);
				\draw (0,0) -- (0,-1) ;

			\end{tikzpicture}
			
			\\
			\hline 
			\begin{tikzpicture}[scale=0.5]
				\path (-1,0) node(x1) {\tiny$\bullet $} ;
				\path (-1,1) node(x2) {\tiny$\bullet$};
				\path (-2,2) node(x3) {\tiny$\bullet$};
				\path (-3,1) node(x4) {\tiny$\bullet$};
				\draw (-1,0) -- (-1,1)node[pos=0.7,right]{\tiny1} ;
				\draw (-1,1) -- (-2,2)node[pos=0.7,right]{\tiny1} ;
				\draw (-2,2) -- (-3,1)node[pos=0.3,left]{\tiny1} ;

			\end{tikzpicture}
			& \begin{tikzpicture}[scale=0.5]
				\path (-2,0) node(x1) {\tiny$\bullet$} ;
				\path (-1,-1) node(x2) {\tiny$\bullet$};
				\path (-2,-2) node(x3) {\tiny$\bullet$};
				
				\draw (-2,0) -- (-1,-1);
				\draw (-1,-1) -- (-2,-2) ;

			\end{tikzpicture}
			& \begin{tikzpicture}[scale=0.5]
				\path (-2,0) node(x1) {\tiny$\bullet$} ;
				\path (-1,-1) node(x2) {\tiny$\bullet$};
				\path (-2,-2) node(x3) {\tiny$\bullet$};
				
				\draw (-2,0) -- (-1,-1);
				\draw (-1,-1) -- (-2,-2) ;

			\end{tikzpicture}
			\\
			\hline
			\begin{tikzpicture}[scale=0.5]
				\path (0,0) node(x1) {\tiny$\bullet$} ;
				\path (-1,-1) node(x2) {\tiny$\bullet$};
				\path (0,-1) node(x3) {\tiny$\bullet$};
				\path (1,-1) node(x4) {\tiny$\bullet$};
				\draw (0,0) -- (-1,-1)node[pos=0.3,left]{\tiny1} ;
				\draw (0,0) -- (0,-1)node[pos=0.6]{\tiny\;\;1} ;
				\draw (0,0) -- (1,-1)node[pos=0.3,right]{\tiny1} ;

			\end{tikzpicture}&
			\begin{tikzpicture}[scale=0.5]
				\path (-3,3) node(x1) {\tiny$\bullet$} ;
				\path (-1,1) node(x2) {\tiny$\bullet$};
				\path (-2,2) node(x3) {\tiny$\bullet$};
				
				\draw (-2,2) -- (-1,1);
				\draw (-3,3) -- (-2,2) ;

			\end{tikzpicture}
			&\begin{tikzpicture}[scale=0.5]
				\path (-3,3) node(x1) {\tiny$\bullet$} ;
				\path (-1,1) node(x2) {\tiny$\bullet$};
				\path (-2,2) node(x3) {\tiny$\bullet$};
				
				\draw (-2,2) -- (-1,1);
				\draw (-3,3) -- (-2,2) ;

			\end{tikzpicture}
			\\
			
			\hline 
			\begin{tikzpicture}[scale=0.5]
				\path (0,0) node(x1) {\tiny$\bullet$} ;
				\path (0,-1) node(x2) {\tiny$\bullet$};
				\path (-1,-2) node(x3) {\tiny$\bullet$};
				\path (1,-2) node(x4) {\tiny$\bullet$};
				\draw (0,0) -- (0,-1)node[pos=0.5,left]{\tiny1} ;
				\draw (0,-1) -- (1,-2)node[pos=0.5,right]{\tiny1} ;
				\draw (0,-1) -- (-1,-2)node[pos=0.5,left]{\tiny 1} ;

			\end{tikzpicture}
			& \begin{tikzpicture}[scale=0.5]
				\path (-2,0) node(x1) {\tiny$\bullet$} ;
				\path (-1,1) node(x2) {\tiny$\bullet$};
				\path (-1,-1) node(x3) {\tiny$\bullet$};
				
				\draw (-2,0) -- (-1,-1);
				\draw (-2,0) -- (-1,1) ;

			\end{tikzpicture}
			& \begin{tikzpicture}[scale=0.5]
				\path (-2,0) node(x1) {\tiny$\bullet$} ;
				\path (-1,1) node(x2) {\tiny$\bullet$};
				\path (-1,-1) node(x3) {\tiny$\bullet$};
				
				\draw (-2,0) -- (-1,-1);
				\draw (-2,0) -- (-1,1) ;

			\end{tikzpicture}\\
			

		\end{tabular}
		
		\caption{First row is a multi-edge tree with 3 edges, second row the corresponding binary tree, according to the classical rotation correspondence, ignoring the unary branches.
			Third row is inserting extra horizontal edges when the multiplicities are higher than 1.}
	\end{table}	
\end{center}


\bibliographystyle{plain}


\end{document}